# On lower limits and equivalences for distribution tails of randomly stopped sums

DENIS DENISOV,[1,2,*] SERGUEI FOSS[2,3,**] and DMITRY KORSHUNOV[3,4,†]

[1]*Eurandom, P.O. Box 513 – 5600 MB Eindhoven, The Netherlands.*
[2]*School of MACS, Heriot-Watt University, Edinburgh EH14 4AS, UK.*
*E-mail:* [*]*Denisov@ma.hw.ac.uk;* [**]*Foss@ma.hw.ac.uk*
[3]*Sobolev Institute of Mathematics, 4 Koptyuga pr., Novosibirsk 630090, Russia.*
*E-mail:* [†]*Korshunov@math.nsc.ru*
[4]*Novosibirsk State University, 2 Pirogova str., Novosibirsk 630090, Russia.*

For a distribution $F^{*\tau}$ of a random sum $S_\tau = \xi_1 + \cdots + \xi_\tau$ of i.i.d. random variables with a common distribution $F$ on the half-line $[0,\infty)$, we study the limits of the ratios of tails $\overline{F^{*\tau}}(x)/\overline{F}(x)$ as $x \to \infty$ (here, $\tau$ is a counting random variable which does not depend on $\{\xi_n\}_{n \geq 1}$). We also consider applications of the results obtained to random walks, compound Poisson distributions, infinitely divisible laws, and subcritical branching processes.

## 1. Introduction

Let $\xi, \xi_1, \xi_2, \ldots$, be independent identically distributed non-negative random variables. We assume that their common distribution $F$ on the half-line $[0, \infty)$ has an unbounded support, that is, $\overline{F}(x) \equiv F(x, \infty) > 0$ for all $x$. Put $S_0 = 0$ and $S_n = \xi_1 + \cdots + \xi_n$, $n = 1, 2, \ldots$.

Let $\tau$ be a counting random variable which does not depend on $\{\xi_n\}_{n \geq 1}$ and which has finite mean. Denote by $F^{*\tau}$ the distribution of a randomly stopped sum $S_\tau = \xi_1 + \cdots + \xi_\tau$.

In this paper, we discuss how the tail behavior of $F^{*\tau}$ relates to that of $F$ and, in particular, under what conditions

$$\liminf_{x \to \infty} \frac{\overline{F^{*\tau}}(x)}{\overline{F}(x)} = \mathbf{E}\tau. \tag{1}$$

Relations on lower limits of ratios of tails were first discussed by Rudin [21]. Theorem $2^*$ of that paper states (for an integer $p$) the following.

**Theorem 1.** *Suppose there exists a positive $p \in [1, \infty)$ such that $\mathbf{E}\xi^p = \infty$, but $\mathbf{E}\tau^p < \infty$. Then* (1) *holds.*

Rudin's studies were motivated by Chover, Ney and Wainger [7] who considered, in particular, the problem of existence of a limit for the ratio

$$\frac{\overline{F^{*\tau}}(x)}{\overline{F}(x)} \qquad \text{as } x \to \infty. \tag{2}$$

From Theorem 1, it follows that if $F$ and $\tau$ satisfy its conditions and if a limit of (2) exists, then that limit must equal $\mathbf{E}\tau$.

Rudin proved Theorem 1 via probability generating function techniques. Below, we give an alternative and more direct proof of Theorem 1 in the case of any positive $p$ (i.e., not necessarily integer). Our method is based on truncation arguments; in this way, we propose a general scheme (see Theorem 4 below) which may also be applied to distributions having all moments finite.

The condition $\mathbf{E}\xi^p = \infty$ rules out many distributions of interest in, say, the theory of subexponential distributions. For example, log-normal and Weibull-type distributions have all moments finite. Our first result presents a natural condition on a stopping time $\tau$ guaranteeing relation (1) for the whole class of heavy-tailed distributions.

Recall that a random variable $\xi$ has a *light-tailed* distribution $F$ on $[0, \infty)$ if $\mathbf{E}e^{\gamma\xi} < \infty$ with some $\gamma > 0$. Otherwise, $F$ is called a *heavy-tailed* distribution; this happens if and only if $\mathbf{E}e^{\gamma\xi} = \infty$ for all $\gamma > 0$.

**Theorem 2.** *Let $F$ be a heavy-tailed distribution and $\tau$ have a light-tailed distribution. Then* (1) *holds.*

The proof of Theorem 2 is based on a new technical tool (see Lemma 2) and significantly differs from the proof of Theorem 1 in Foss and Korshunov [15], where the particular case $\tau = 2$ was considered. Theorem 2 is restricted to the case of light-tailed $\tau$, but here, we extend Rudin's result to the class of all heavy-tailed distributions. The reasons for the restriction to $\mathbf{E}e^{\gamma\tau} < \infty$ come from the proof of Theorem 2, but are, in fact, rather natural: the tail of $\tau$ should be lighter than the tail of any heavy-tailed distribution. Indeed, if $\xi_1 \geq 1$, then $\overline{F^{*\tau}}(x) \geq \mathbf{P}\{\tau > x\}$. This shows that the tail of $F^{*\tau}$ is at least as heavy as that of $\tau$. Note that in Theorem 1, in some sense, the tail of $F$ is heavier than the tail of $\tau$.

Theorem 2 may be applied in various areas where randomly stopped sums appear; see Sections 8–11 (random walks, compound Poisson distributions, infinitely divisible laws and branching processes) and, for instance, Kalashnikov [17] for further examples.

For any distribution on $[0, \infty)$, let

$$\varphi(\gamma) = \int_0^\infty e^{\gamma x} F(\mathrm{d}x) \in (0, \infty], \qquad \gamma \in \mathbf{R},$$

and
$$\widehat{\gamma} = \sup\{\gamma : \varphi(\gamma) < \infty\} \in [0, \infty].$$

Note that the moment generating function $\varphi(\gamma)$ is increasing and continuous in the interval $(-\infty, \widehat{\gamma})$ and that $\varphi(\widehat{\gamma}) = \lim_{\gamma \uparrow \widehat{\gamma}} \varphi(\gamma) \in [1, \infty]$. The following result was proven in Foss and Korshunov [15], Theorem 3. Let
$$\frac{\overline{F * F}(x)}{\overline{F}(x)} \to c \quad \text{as } x \to \infty,$$

where $c \in (0, \infty]$. Then, necessarily, $c = 2\varphi(\widehat{\gamma})$. We state now a generalization to $\tau$-fold convolution.

**Theorem 3.** *Let $\varphi(\widehat{\gamma}) < \infty$ and $\mathbf{E}(\varphi(\widehat{\gamma}) + \varepsilon)^\tau < \infty$ for some $\varepsilon > 0$. Assume that*
$$\frac{\overline{F^{*\tau}}(x)}{\overline{F}(x)} \to c \quad \text{as } x \to \infty,$$

*where $c \in (0, \infty]$. Then $c = \mathbf{E}(\tau \varphi^{\tau-1}(\widehat{\gamma}))$.*

For (comments on) earlier partial results in the case $\tau = 2$, see, for example, Chover, Ney and Wainger [6, 7], Cline [8], Embrechts and Goldie [10], Foss and Korshunov [15], Pakes [19], Rogozin [20], Teugels [23] and further references therein. The proof of Theorem 3 follows from Lemmas 3 and 4 in Section 7.

## 2. Preliminary result

We start with the following result.

**Theorem 4.** *Assume that there exists a non-decreasing concave function $h : \mathbf{R}^+ \to \mathbf{R}^+$ such that*
$$\mathbf{E}\mathrm{e}^{h(\xi)} < \infty \quad \text{and} \quad \mathbf{E}\xi \mathrm{e}^{h(\xi)} = \infty. \tag{3}$$

*For any $n \geq 1$, put $A_n = \mathbf{E}\mathrm{e}^{h(\xi_1 + \cdots + \xi_n)}$. Assume that $F$ is heavy-tailed and that*
$$\mathbf{E}\tau A_{\tau-1} < \infty. \tag{4}$$

*Then, for any light-tailed distribution $G$ on $[0, \infty)$,*
$$\liminf_{x \to \infty} \frac{\overline{G * F^{*\tau}}(x)}{\overline{F}(x)} = \mathbf{E}\tau. \tag{5}$$

By considering $G$ concentrated at 0, we get the following.

**Corollary 1.** *In the conditions of Theorem 4, (1) holds.*

In order to prove Theorem 4, first we restate Theorem 1$^*$ of Rudin [21] (in Lemma 1 below) in terms of probability distributions and stopping times.

**Lemma 1.** *For any distribution $F$ on $[0, \infty)$ with unbounded support and any counting random variable $\tau$,*

$$\liminf_{x \to \infty} \frac{\overline{F^{*\tau}}(x)}{\overline{F}(x)} \geq \mathbf{E}\tau.$$

**Proof.** For any two distributions $F_1$ and $F_2$ on $[0, \infty)$ with unbounded supports,

$$\overline{F_1 * F_2}(x) \geq (F_1 \times F_2)((x, \infty) \times [0, x]) + (F_1 \times F_2)([0, x] \times (x, \infty))$$
$$\sim \overline{F}_1(x) + \overline{F}_2(x) \qquad \text{as } x \to \infty.$$

By induction arguments, this implies that, for any $n \geq 1$,

$$\liminf_{x \to \infty} \frac{\overline{F^{*n}}(x)}{\overline{F}(x)} \geq n.$$

Applying Fatou's lemma to the representation

$$\frac{\overline{F^{*\tau}}(x)}{\overline{F}(x)} = \sum_{n=1}^{\infty} \mathbf{P}\{\tau = n\} \frac{\overline{F^{*n}}(x)}{\overline{F}(x)},$$

completes the proof. □

**Proof of Theorem 4.** It follows from Lemma 1 that it is sufficient to prove the following inequality:

$$\liminf_{x \to \infty} \frac{\overline{G * F^{*\tau}}(x)}{\overline{F}(x)} \leq \mathbf{E}\tau.$$

Assume the contrary, that is, that there exist $\delta > 0$ and $x_0$ such that

$$\overline{G * F^{*\tau}}(x) \geq (\mathbf{E}\tau + \delta)\overline{F}(x) \qquad \text{for all } x > x_0. \tag{6}$$

For any positive $b > 0$, consider a concave function

$$h_b(x) \equiv \min\{h(x), bx\}, \tag{7}$$

which is non-negative because $h \geq 0$. Since $F$ is heavy-tailed, $h(x) = o(x)$ as $x \to \infty$. Therefore, for any fixed $b$, there exists $x_1$ such that $h_b(x) = h(x)$ for all $x > x_1$. Hence, by condition (3),

$$\mathbf{E}e^{h_b(\xi)} < \infty \quad \text{and} \quad \mathbf{E}\xi e^{h_b(\xi)} = \infty. \tag{8}$$

For any $x$, we have the convergence $h_b(x) \downarrow 0$ as $b \downarrow 0$. Then, for any fixed $n$,
$$A_{n,b} \equiv \mathbf{E} e^{h_b(\xi_1 + \cdots + \xi_n)} \downarrow 1 \quad \text{as } b \downarrow 0.$$
This and condition (4) together imply that there exists $b$ such that
$$\mathbf{E}\tau A_{\tau-1,b} \leq \mathbf{E}\tau + \delta/8. \tag{9}$$

Let $\eta$ be a random variable with distribution $G$ which does not depend on $\{\xi_n\}_{n \geq 1}$ and $\tau$. Since $G$ is light-tailed,
$$\mathbf{E}\eta e^{h_b(\eta)} < \infty. \tag{10}$$
In addition, we may choose $b > 0$ sufficiently small that
$$\mathbf{E} e^{h_b(\eta)} (\mathbf{E}\tau + \delta/8) \leq \mathbf{E}\tau + \delta/4. \tag{11}$$

For any real $a$ and $t$, put $a^{[t]} = \min\{a, t\}$. Then
$$\frac{\mathbf{E}(\eta + \xi_1^{[t]} + \cdots + \xi_\tau^{[t]}) e^{h_b(\eta + \xi_1 + \cdots + \xi_\tau)}}{\mathbf{E}\xi_1^{[t]} e^{h_b(\xi_1)}} = \sum_{n=1}^{\infty} \frac{\mathbf{E}\eta e^{h_b(\eta + \xi_1 + \cdots + \xi_n)}}{\mathbf{E}\xi_1^{[t]} e^{h_b(\xi_1)}} \mathbf{P}\{\tau = n\}$$
$$+ \sum_{n=1}^{\infty} n \frac{\mathbf{E}\xi_1^{[t]} e^{h_b(\eta + \xi_1 + \cdots + \xi_n)}}{\mathbf{E}\xi_1^{[t]} e^{h_b(\xi_1)}} \mathbf{P}\{\tau = n\}.$$

By the concavity of the function $h_b$,
$$\sum_{n=1}^{\infty} \frac{\mathbf{E}\eta e^{h_b(\eta + \xi_1 + \cdots + \xi_n)}}{\mathbf{E}\xi_1^{[t]} e^{h_b(\xi_1)}} \mathbf{P}\{\tau = n\} \leq \sum_{n=1}^{\infty} \frac{\mathbf{E}\eta e^{h_b(\eta) + h_b(\xi_1 + \cdots + \xi_n)}}{\mathbf{E}\xi_1^{[t]} e^{h_b(\xi_1)}} \mathbf{P}\{\tau = n\}$$
$$= \frac{\mathbf{E}\eta e^{h_b(\eta)}}{\mathbf{E}\xi_1^{[t]} e^{h_b(\xi_1)}} \mathbf{E} A_{\tau,b}$$
$$\to 0 \quad \text{as } t \to \infty,$$
due to (10), (9) and (8). Again, by the concavity of the function $h_b$,
$$\sum_{n=1}^{\infty} n \frac{\mathbf{E}\xi_1^{[t]} e^{h_b(\eta + \xi_1 + \cdots + \xi_n)}}{\mathbf{E}\xi_1^{[t]} e^{h_b(\xi_1)}} \mathbf{P}\{\tau = n\} \leq \sum_{n=1}^{\infty} n \frac{\mathbf{E}\xi_1^{[t]} e^{h_b(\eta) + h_b(\xi_1) + h_b(\xi_2 + \cdots + \xi_n)}}{\mathbf{E}\xi_1^{[t]} e^{h_b(\xi_1)}} \mathbf{P}\{\tau = n\}$$
$$= \mathbf{E} e^{h_b(\eta)} \sum_{n=1}^{\infty} n A_{n-1,b} \mathbf{P}\{\tau = n\}$$
$$\leq \mathbf{E}\tau + \delta/4,$$
by (9) and (11). Hence, for sufficiently large $t$,
$$\frac{\mathbf{E}(\eta + \xi_1^{[t]} + \cdots + \xi_\tau^{[t]}) e^{h_b(\eta + \xi_1 + \cdots + \xi_\tau)}}{\mathbf{E}\xi_1^{[t]} e^{h_b(\xi_1)}} \leq \mathbf{E}\tau + \delta/2. \tag{12}$$

On the other hand, since $(\eta + \xi_1 + \cdots + \xi_\tau)^{[t]} \leq \eta + \xi_1^{[t]} + \cdots + \xi_\tau^{[t]}$,

$$\frac{\mathbf{E}(\eta + \xi_1^{[t]} + \cdots + \xi_\tau^{[t]})e^{h_b(\eta+\xi_1+\cdots+\xi_\tau)}}{\mathbf{E}\xi_1^{[t]}e^{h_b(\xi_1)}} \geq \frac{\mathbf{E}(\eta + \xi_1 + \cdots + \xi_\tau)^{[t]}e^{h_b(\eta+\xi_1+\cdots+\xi_\tau)}}{\mathbf{E}\xi_1^{[t]}e^{h_b(\xi_1)}}$$

$$= \frac{\int_0^\infty x^{[t]}e^{h_b(x)}(G * F^{*\tau})(\mathrm{d}x)}{\int_0^\infty x^{[t]}e^{h_b(x)}F(\mathrm{d}x)}. \qquad (13)$$

The right-hand side, after integration by parts, is equal to

$$\frac{\int_0^\infty \overline{G * F^{*\tau}}(x)\,\mathrm{d}(x^{[t]}e^{h_b(x)})}{\int_0^\infty \overline{F}(x)\,\mathrm{d}(x^{[t]}e^{h_b(x)})}.$$

Since $\mathbf{E}\xi_1 e^{h_b(\xi_1)} = \infty$, both integrals in this fraction tend to infinity as $t \to \infty$. For the *non-decreasing* function $h_b(x)$, the latter fact and assumption (6) together imply that

$$\liminf_{t\to\infty} \frac{\int_0^\infty \overline{G * F^{*\tau}}(x)\,\mathrm{d}(x^{[t]}e^{h_b(x)})}{\int_0^\infty \overline{F}(x)\,\mathrm{d}(x^{[t]}e^{h_b(x)})} = \liminf_{t\to\infty} \frac{\int_{x_0}^\infty \overline{G * F^{*\tau}}(x)\,\mathrm{d}(x^{[t]}e^{h_b(x)})}{\int_{x_0}^\infty \overline{F}(x)\,\mathrm{d}(x^{[t]}e^{h_b(x)})} \geq \mathbf{E}\tau + \delta.$$

Substituting this into (13), we get a contradiction of (12) for sufficiently large $t$. The proof is thus complete. $\square$

## 3. Proof of Theorem 1

Take an integer $k \geq 0$ such that $p - 1 \leq k < p$. Without loss of generality, we may assume that $\mathbf{E}\xi^k < \infty$ (otherwise, we may consider a smaller $p$).

Consider a concave non-decreasing function $h(x) = (p-1)\ln x$. Then $\mathbf{E}e^{h(\xi_1)} < \infty$ and $\mathbf{E}\xi_1 e^{h(\xi_1)} = \infty$. Thus,

$$A_n \equiv \mathbf{E}e^{h(\xi_1+\cdots+\xi_n)} = \mathbf{E}(\xi_1 + \cdots + \xi_n)^{p-1}$$
$$\leq (\mathbf{E}(\xi_1 + \cdots + \xi_n)^k)^{(p-1)/k}$$

since $(p-1)/k \leq 1$. Further,

$$\mathbf{E}(\xi_1 + \cdots + \xi_n)^k = \sum_{i_1,\ldots,i_k=1}^n \mathbf{E}(\xi_{i_1} \cdots \xi_{i_k}) \leq cn^k,$$

where

$$c \equiv \sup_{1 \leq i_1,\ldots,i_k \leq n} \mathbf{E}(\xi_{i_1} \cdots \xi_{i_k}) < \infty,$$

due to the fact that $\mathbf{E}\xi^k < \infty$. Hence, $A_n \leq c^{(p-1)/k} n^{p-1}$ for all $n$. Therefore, we get $\mathbf{E}\tau A_{\tau-1} \leq c^{(p-1)/k}\mathbf{E}\tau^p < \infty$. All conditions of Theorem 4 are met and the proof is complete.

## 4. Characterization of heavy-tailed distributions

In the sequel, we need the following existence result which strengthens a lemma in Rudin [21], page 989; and Lemma 1 in Foss and Korshunov [15]. Fix any $\delta \in (0,1]$.

**Lemma 2.** *If a random variable $\xi \geq 0$ has a heavy-tailed distribution, then there exists a non-decreasing concave function $h : \mathbf{R}^+ \to \mathbf{R}^+$ such that $\mathbf{E} e^{h(\xi)} \leq 1 + \delta$ and $\mathbf{E}\xi e^{h(\xi)} = \infty$.*

**Proof.** Without loss of generality, assume that $\xi > 0$ a.s., that is, that $\overline{F}(0) = 1$. We will construct a piecewise linear function $h(x)$. For that, we introduce two positive sequences, $x_n \uparrow \infty$ and $\varepsilon_n \downarrow 0$ as $n \to \infty$, and let

$$h(x) = h(x_{n-1}) + \varepsilon_n (x - x_{n-1}) \qquad \text{if } x \in (x_{n-1}, x_n], \ n \geq 1.$$

This function is non-decreasing since $\varepsilon_n > 0$. Moreover, this function is concave due to the monotonicity of $\varepsilon_n$.

Put $x_0 = 0$ and $h(0) = 0$. Since $\xi$ is heavy-tailed, we can choose $x_1 \geq 2$ so that

$$\mathbf{E}\{e^{\xi}; \xi \in (x_0, x_1]\} + e^{x_1} \overline{F}(x_1) > 1 + \delta.$$

Choose $\varepsilon_1 > 0$ so that

$$\mathbf{E}\{e^{\varepsilon_1 \xi}; \xi \in (x_0, x_1]\} + e^{\varepsilon_1 x_1} \overline{F}(x_1) = e^{h(x_0)} \overline{F}(0) + \delta/2 = 1 + \delta/2,$$

which is equivalent to

$$\mathbf{E}\{e^{h(\xi)}; \xi \in (x_0, x_1]\} + e^{h(x_1)} \overline{F}(x_1) = e^{h(x_0)} \overline{F}(0) + \delta/2.$$

By induction, we construct an increasing sequence $x_n$ and a decreasing sequence $\varepsilon_n > 0$ such that $x_n \geq 2^n$ and

$$\mathbf{E}\{e^{h(\xi)}; \xi \in (x_{n-1}, x_n]\} + e^{h(x_n)} \overline{F}(x_n) = e^{h(x_{n-1})} \overline{F}(x_{n-1}) + \delta/2^n$$

for any $n \geq 2$. For $n = 1$, this is already done. Make the induction hypothesis for some $n \geq 2$. Due to heavy-tailedness, there exists $x_{n+1} \geq 2^{n+1}$ sufficiently large that

$$\mathbf{E}\{e^{\varepsilon_n (\xi - x_n)}; \xi \in (x_n, x_{n+1}]\} + e^{\varepsilon_n (x_{n+1} - x_n)} \overline{F}(x_{n+1}) > 1 + \delta.$$

Note that

$$\mathbf{E}\{e^{\varepsilon_{n+1}(\xi - x_n)}; \xi \in (x_n, x_{n+1}]\} + e^{\varepsilon_{n+1}(x_{n+1} - x_n)} \overline{F}(x_{n+1})$$

as a function of $\varepsilon_{n+1}$ is continuously decreasing to $\overline{F}(x_n)$ as $\varepsilon_{n+1} \downarrow 0$. Therefore, we can choose $\varepsilon_{n+1} \in (0, \varepsilon_n)$ so that

$$\mathbf{E}\{e^{\varepsilon_{n+1}(\xi - x_n)}; \xi \in (x_n, x_{n+1}]\} + e^{\varepsilon_{n+1}(x_{n+1} - x_n)} \overline{F}(x_{n+1})$$
$$= \overline{F}(x_n) + \delta/(2^{n+1} e^{h(x_n)}).$$

By definition of $h(x)$, this is equivalent to the following equality:

$$\mathbf{E}\{e^{h(\xi)}; \xi \in (x_n, x_{n+1}]\} + e^{h(x_{n+1})}\overline{F}(x_{n+1}) = e^{h(x_n)}\overline{F}(x_n) + \delta/2^{n+1}.$$

Our induction hypothesis now holds with $n+1$ in place of $n$, as required.

Next,

$$\mathbf{E}e^{h(\xi)} = \sum_{n=0}^{\infty} \mathbf{E}\{e^{h(\xi)}; \xi \in (x_n, x_{n+1}]\}$$

$$= \sum_{n=0}^{\infty}(e^{h(x_n)}\overline{F}(x_n) - e^{h(x_{n+1})}\overline{F}(x_{n+1}) + \delta/2^{n+1})$$

$$= e^{h(x_0)}\overline{F}(x_0) + \delta = 1 + \delta.$$

On the other hand, since $x_k \geq 2^k$,

$$\mathbf{E}\{\xi e^{h(\xi)}; \xi > x_n\} = \sum_{k=n}^{\infty} \mathbf{E}\{\xi e^{h(\xi)}; \xi \in (x_k, x_{k+1}]\}$$

$$\geq 2^n \sum_{k=n}^{\infty} \mathbf{E}\{e^{h(\xi)}; \xi \in (x_k, x_{k+1}]\}$$

$$\geq 2^n \sum_{k=n}^{\infty}(e^{h(x_k)}\overline{F}(x_k) - e^{h(x_{k+1})}\overline{F}(x_{k+1}) + \delta/2^{k+1}).$$

Then, for any $n$,

$$\mathbf{E}\{\xi e^{h(\xi)}; \xi > x_n\} \geq 2^n(e^{h(x_n)}\overline{F}(x_n) + \delta/2^n) \geq \delta,$$

which implies that $\mathbf{E}\xi e^{h(\xi)} = \infty$. Also note that, necessarily, $\lim_{n \to \infty} \varepsilon_n = 0$; otherwise, $\liminf_{x \to \infty} h(x)/x > 0$ and $\xi$ is light-tailed. The proof of the lemma is thus complete. □

## 5. Proof of Theorem 2

Since $\tau$ has a light-tailed distribution,

$$\mathbf{E}\tau(1+\varepsilon)^{\tau-1} < \infty$$

for some sufficiently small $\varepsilon > 0$. By Lemma 2, there exists a concave increasing function $h$, $h(0) = 0$, such that $\mathbf{E}e^{h(\xi_1)} \leq 1 + \varepsilon$ and $\mathbf{E}\xi_1 e^{h(\xi_1)} = \infty$. Then, by concavity,

$$A_n \equiv \mathbf{E}e^{h(\xi_1 + \cdots + \xi_n)} \leq \mathbf{E}e^{h(\xi_1) + \cdots + h(\xi_n)} \leq (1+\varepsilon)^n.$$

Combining, we get $\mathbf{E}\tau A_{\tau-1} < \infty$. All conditions of Theorem 4 are met and the proof is thus complete.

## 6. Fractional exponential moments

One can go further and obtain various results on lower limits and equivalences for heavy-tailed distributions $F$ which have all finite power moments (e.g., Weibull and log-normal distributions). For instance, we have the following result (see Denisov, Foss and Korshunov [9] for the proof).

Suppose there exists $\alpha$, $0 < \alpha < 1$, such that $\mathbf{E}e^{c\xi^\alpha} = \infty$ for all $c > 0$. If $\mathbf{E}e^{\delta\tau^\alpha} < \infty$ for some $\delta > 0$, then (1) holds.

## 7. Tail equivalence for randomly stopped sums

The following auxiliary lemma compares the tail behavior of the convolution tail and that of the exponentially transformed distribution.

**Lemma 3.** *Let the distribution $F$ and the number $\gamma \geq 0$ be such that $\varphi(\gamma) < \infty$. Let the distribution $G$ be the result of the exponential change of measure with parameter $\gamma$, that is, $G(\mathrm{d}u) = e^{\gamma u} F(\mathrm{d}u)/\varphi(\gamma)$. Let $\tau$ be any counting random variable such that $\mathbf{E}\varphi^\tau(\gamma) < \infty$ and let $\nu$ have the distribution $\mathbf{P}\{\nu = k\} = \varphi^k(\gamma)\mathbf{P}\{\tau = k\}/\mathbf{E}\varphi^\tau(\gamma)$. Then*

$$\liminf_{x \to \infty} \frac{\overline{G^{*\nu}}(x)}{\overline{G}(x)} \geq \frac{1}{\mathbf{E}\varphi^{\tau-1}(\gamma)} \liminf_{x \to \infty} \frac{\overline{F^{*\tau}}(x)}{\overline{F}(x)}$$

*and*

$$\limsup_{x \to \infty} \frac{\overline{G^{*\nu}}(x)}{\overline{G}(x)} \leq \frac{1}{\mathbf{E}\varphi^{\tau-1}(\gamma)} \limsup_{x \to \infty} \frac{\overline{F^{*\tau}}(x)}{\overline{F}(x)}.$$

**Proof.** Put

$$\widehat{c} \equiv \liminf_{x \to \infty} \frac{\overline{F^{*\tau}}(x)}{\overline{F}(x)}.$$

By Lemma 1, $\widehat{c} \in [\mathbf{E}\tau, \infty]$. For any fixed $c \in (0, \widehat{c})$, there exists $x_0 > 0$ such that, for any $x > x_0$,

$$\overline{F^{*\tau}}(x) \geq c\overline{F}(x). \tag{14}$$

By the total probability law,

$$\overline{G^{*\nu}}(x) = \sum_{k=1}^{\infty} \mathbf{P}\{\nu = k\}\overline{G^{*k}}(x)$$

$$= \sum_{k=1}^{\infty} \frac{\varphi^k(\gamma)\mathbf{P}\{\tau = k\}}{\mathbf{E}\varphi^\tau(\gamma)} \int_x^\infty e^{\gamma y} \frac{F^{*k}(\mathrm{d}y)}{\varphi^k(\gamma)}$$

$$= \frac{1}{\mathbf{E}\varphi^\tau(\gamma)} \sum_{k=1}^{\infty} \mathbf{P}\{\tau = k\} \int_x^\infty e^{\gamma y} F^{*k}(\mathrm{d}y).$$

Integrating by parts, we obtain

$$\sum_{k=1}^{\infty} \mathbf{P}\{\tau=k\} \left[ e^{\gamma x} \overline{F^{*k}}(x) + \int_x^\infty \overline{F^{*k}}(y) \, de^{\gamma y} \right] = e^{\gamma x} \overline{F^{*\tau}}(x) + \int_x^\infty \overline{F^{*\tau}}(y) \, de^{\gamma y}.$$

Also using (14) we get, for $x > x_0$,

$$\overline{G^{*\nu}}(x) \geq \frac{c}{\mathbf{E}\varphi^\tau(\gamma)} \left[ e^{\gamma x} \overline{F}(x) + \int_x^\infty \overline{F}(y) \, de^{\gamma y} \right]$$

$$= \frac{c}{\mathbf{E}\varphi^\tau(\gamma)} \int_x^\infty e^{\gamma y} F(dy) = \frac{c}{\mathbf{E}\varphi^{\tau-1}(\gamma)} \overline{G}(x).$$

Letting $c \uparrow \widehat{c}$, we obtain the first conclusion of the lemma. The proof of the second conclusion follows similarly. □

**Lemma 4.** *If $0 < \widehat{\gamma} < \infty$, $\varphi(\widehat{\gamma}) < \infty$ and $\mathbf{E}(\varphi(\widehat{\gamma}) + \varepsilon)^\tau < \infty$ for some $\varepsilon > 0$, then*

$$\liminf_{x \to \infty} \frac{\overline{F^{*\tau}}(x)}{\overline{F}(x)} \leq \mathbf{E}\tau \varphi^{\tau-1}(\widehat{\gamma})$$

*and*

$$\limsup_{x \to \infty} \frac{\overline{F^{*\tau}}(x)}{\overline{F}(x)} \geq \mathbf{E}\tau \varphi^{\tau-1}(\widehat{\gamma}).$$

**Proof.** We apply the exponential change of measure with parameter $\widehat{\gamma}$ and consider the distribution $G(du) = e^{\widehat{\gamma} u} F(du)/\varphi(\widehat{\gamma})$ and the stopping time $\nu$ with the distribution $\mathbf{P}\{\nu = k\} = \varphi^k(\widehat{\gamma}) \mathbf{P}\{\tau = k\}/\mathbf{E}\varphi^\tau(\widehat{\gamma})$. From the definition of $\widehat{\gamma}$, the distribution $G$ is heavy-tailed. The distribution of $\nu$ is light-tailed because $\mathbf{E}e^{\kappa \nu} < \infty$ with $\kappa = \ln(\varphi(\widehat{\gamma}) + \varepsilon) - \ln \varphi(\widehat{\gamma}) > 0$. Hence,

$$\limsup_{x \to \infty} \frac{\overline{G^{*\nu}}(x)}{\overline{G}(x)} \geq \liminf_{x \to \infty} \frac{\overline{G^{*\nu}}(x)}{\overline{G}(x)} = \mathbf{E}\nu,$$

by Theorem 2. The result now follows from Lemma 3 with $\gamma = \widehat{\gamma}$, since $\mathbf{E}\nu = \mathbf{E}\tau \varphi^\tau(\widehat{\gamma})/\mathbf{E}\varphi^\tau(\widehat{\gamma})$. □

**Proof of Theorem 3.** In the case where $F$ is heavy-tailed, we have $\widehat{\gamma} = 0$ and $\varphi(\widehat{\gamma}) = 1$. By Theorem 2, $c = \mathbf{E}\tau$, as required.

In the case $\widehat{\gamma} \in (0, \infty)$ and $\varphi(\widehat{\gamma}) < \infty$, the desired conclusion follows from Lemma 4. □

## 8. Supremum of a random walk

Hereafter, we need the notion of subexponential distributions. A distribution $F$ on $\mathbf{R}^+$ is called *subexponential* if $\overline{F * F}(x) \sim 2\overline{F}(x)$ as $x \to \infty$.

Let $\{\xi_n\}$ be a sequence of independent random variables with a common distribution $F$ on $\mathbf{R}$ and $\mathbf{E}\xi_1 = -m < 0$. Put $S_0 = 0$, $S_n = \xi_1 + \cdots + \xi_n$. By the strong law of large numbers (SLLN), $M = \sup_{n \geq 0} S_n$ is finite with probability 1.

Let $F^I$ be the integrated tail distribution on $\mathbf{R}^+$, that is,

$$\overline{F^I}(x) \equiv \min\left(1, \int_x^\infty \overline{F}(y)\,\mathrm{d}y\right), \qquad x > 0.$$

It is well known (see, e.g., Asmussen [1], Embrechts, Klüppelberg and Mikosch [12], Embrechts and Veraverbeke [13] and references therein) that if $F^I$ is subexponential, then

$$\mathbf{P}\{M > x\} \sim \frac{1}{m}\overline{F^I}(x) \qquad \text{as } x \to \infty. \tag{15}$$

Korshunov [18] proved the converse: (15) implies subexponentiality of $F^I$. We now supplement this assertion with the following result.

**Theorem 5.** *Let $F^I$ be long-tailed, that is, $\overline{F^I}(x+1) \sim \overline{F^I}(x)$ as $x \to \infty$. If, for some $c > 0$,*

$$\mathbf{P}\{M > x\} \sim c\overline{F^I}(x) \qquad \text{as } x \to \infty,$$

*then $c = 1/m$ and $F^I$ is subexponential.*

**Proof.** Consider the defective stopping time

$$\eta = \inf\{n \geq 1 : S_n > 0\} \leq \infty$$

and let $\{\psi_n\}$ be i.i.d. random variables with common distribution function

$$G(x) \equiv \mathbf{P}\{\psi_n \leq x\} = \mathbf{P}\{S_\eta \leq x \mid \eta < \infty\}.$$

It is well known (see, e.g., Feller [14], Chapter XII) that the distribution of the maximum $M$ coincides with the distribution of the randomly stopped sum $\psi_1 + \cdots + \psi_\tau$, where the counting random variable $\tau$ is independent of the sequence $\{\psi_n\}$ and is geometrically distributed with parameter $p = \mathbf{P}\{M > 0\} < 1$, that is, $\mathbf{P}\{\tau = k\} = (1-p)p^k$ for $k = 0, 1, \ldots$. Equivalently,

$$\mathbf{P}\{M \in B\} = G^{*\tau}(B).$$

It follows from Borovkov [4], Chapter 4, Theorem 10, that if $F^I$ is long-tailed, then

$$\overline{G}(x) \sim \frac{1-p}{pm}\overline{F^I}(x). \tag{16}$$

The theorem hypothesis then implies that

$$\overline{G^{*\tau}}(x) \sim \frac{cpm}{1-p}\overline{G}(x) \qquad \text{as } x \to \infty.$$

Therefore, by Theorem 3 with $\widehat{\gamma} = 0$, $c = \mathbf{E}\tau(1-p)/pm = 1/m$. It then follows from Korshunov [18] that $F^I$ is subexponential. The proof is now complete. □

## 9. The compound Poisson distribution

Let $F$ be a distribution on $\mathbf{R}_+$ and $t$ a positive constant. Let $G$ be the compound Poisson distribution

$$G = \mathrm{e}^{-t} \sum_{n \geq 0} \frac{t^n}{n!} F^{*n}.$$

Considering $\tau$ in Theorem 3 with distribution $\mathbf{P}\{\tau = n\} = t^n \mathrm{e}^{-t}/n!$, we get the following result.

**Theorem 6.** *Let $\varphi(\widehat{\gamma}) < \infty$. If, for some $c > 0$, $\overline{G}(x) \sim c\overline{F}(x)$ as $x \to \infty$, then $c = t\mathrm{e}^{t(\varphi(\widehat{\gamma})-1)}$.*

**Corollary 2.** *The following statements are equivalent:*

   (i) *$F$ is subexponential;*
  (ii) *$G$ is subexponential;*
 (iii) *$\overline{G}(x) \sim t\overline{F}(x)$ as $x \to \infty$;*
 (iv) *$F$ is heavy-tailed and $\overline{G}(x) \sim c\overline{F}(x)$ as $x \to \infty$, for some $c > 0$.*

**Proof.** Equivalence of (i), (ii) and (iii) was proven in Embrechts, Goldie and Veraverbeke [11], Theorem 3. The implication (iv) ⇒ (iii) follows from Theorem 3 with $\widehat{\gamma} = 0$. □

Some local aspects of this problem for heavy-tailed distributions were discussed in Asmussen, Foss and Korshunov [2], Theorem 6.

## 10. Infinitely divisible laws

Let $H$ be an infinitely divisible law on $[0, \infty)$. The Laplace transform of an infinitely divisible law $F$ can be expressed as

$$\int_0^\infty \mathrm{e}^{-\lambda x} H(\mathrm{d}x) = \mathrm{e}^{-a\lambda - \int_0^\infty (1-\mathrm{e}^{-\lambda x})\nu(\mathrm{d}x)}$$

(see, e.g., Feller [14], Chapter XVII). Here, $a \geq 0$ is a constant and the Lévy measure $\nu$ is a Borel measure on $(0, \infty)$ with the properties $\mu = \nu(1, \infty) < \infty$ and $\int_0^1 x\nu(\mathrm{d}x) < \infty$. Put $F(B) = \nu(B \cap (1, \infty))/\mu$.

Relations between the tail behavior of measure $H$ and of the corresponding Lévy measure $\nu$ were considered in Embrechts, Goldie and Veraverbeke [11], Pakes [19] and Shimura and Watanabe [22]. The local analog of that result was proven in Asmussen,

Foss and Korshunov [2]. We strengthen the corresponding result of Embrechts, Goldie and Veraverbeke [11] in the following way.

**Theorem 7.** *The following assertions are equivalent:*

(i) $H$ *is subexponential;*
(ii) $F$ *is subexponential;*
(iii) $\overline{\nu}(x) \sim \overline{H}(x)$ *as* $x \to \infty$;
(iv) $H$ *is heavy-tailed and* $\overline{\nu}(x) \sim c\overline{H}(x)$ *as* $x \to \infty$, *for some* $c > 0$.

**Proof.** Equivalence of (i), (ii) and (iii) was proven in Embrechts, Goldie and Veraverbeke [11], Theorem 1.

It remains to prove the implication (iv) $\Rightarrow$ (iii). It is pointed out in Embrechts, Goldie and Veraverbeke [11] that the distribution $H$ admits the representation $H = G * F^{*\tau}$, where $\overline{G}(x) = O(\mathrm{e}^{-\varepsilon x})$ for some $\varepsilon > 0$ and $\tau$ has a Poisson distribution with parameter $\mu$. Since $H$ is heavy-tailed and $G$ is light-tailed, $F$ is necessarily heavy-tailed. Then, by Theorem 4, we get

$$\liminf_{x \to \infty} \frac{\overline{H}(x)}{\overline{F}(x)} \equiv \liminf_{x \to \infty} \frac{\overline{G * F^{*\tau}}(x)}{\overline{F}(x)} = \mathbf{E}\tau = \mu.$$

On the other hand, for $x > 1$,

$$\frac{\overline{H}(x)}{\overline{F}(x)} = \mu \frac{\overline{H}(x)}{\overline{\nu}(x)} \to \mu c \qquad \text{as } x \to \infty,$$

by assumption (iv). Hence, $c = 1$. □

## 11. Branching processes

In this section, we consider the limit behavior of subcritical, age-dependent branching processes for which the Malthusian parameter does not exist.

Let $h(z)$ be the particle production generating function of an age-dependent branching process with particle lifetime distribution $F$ (see Athreya and Ney [3], Chapter IV, Harris [16], Chapter VI for background). We take the process to be subcritical, that is, $A \equiv h'(1) < 1$. Let $Z(t)$ denote the number of particles at time $t$. It is known (see, e.g., Athreya and Ney [3], Chapter IV, Section 5, or Chistyakov [5]) that $\mathbf{E}Z(t)$ admits the representation

$$\mathbf{E}Z(t) = (1 - A) \sum_{n=1}^{\infty} A^{n-1} \overline{F^{*n}}(t).$$

It was proven in Chistyakov [5] for sufficiently small values of $A$ and then in Chover, Ney and Wainger [6, 7] for any $A < 1$ that $\mathbf{E}Z(t) \sim \overline{F}(t)/(1 - A)$ as $t \to \infty$, provided $F$ is subexponential. The local asymptotics were considered in Asmussen, Foss and Korshunov [2].

Applying Theorem 3 with $\tau$ geometrically distributed and $\widehat{\gamma} = 0$, we deduce the following.

**Theorem 8.** *Let $F$ be heavy-tailed, and, for some $c > 0$, $\mathbf{E}Z(t) \sim c\overline{F}(t)$ as $t \to \infty$. Then $c = 1/(1-A)$ and $F$ is subexponential.*

## Acknowledgements

We are grateful to Toshiro Watanabe who pointed out some incorrectness in the first version of the proof of Theorem 7.

Denisov was supported by the Dutch BSIK project (*BRICKS*). The research of Denisov and Foss was partially supported by EURO-NGI Framework 6 grant. The research of Foss and Korshunov was partially supported by the Royal Society International Joint Project Grant 2005/R2 UP. Korshunov wishes to thank the Russian Science Support Foundation.